\documentclass[12pt,twoside,leqno]{article}

\author{\footnotesize D.D.\ PORO\c SNIUC}
\date{}
\title{\normalsize \bf A CLASS OF K\"AHLER EINSTEIN STRUCTURES ON THE NONZERO
COTANGENT BUNDLE OF A SPACE FORM}

\begin{document}

\maketitle \pagestyle{myheadings}\markboth{\small D.D.Poro\c sniuc
\hspace{5.2cm}}{\hspace{3.3cm} \small A class of K\"ahler Einstein
structures}

\begin{flushleft}

{\scriptsize We obtain a class of K\"ahler Einstein structures on
the nonzero cotangent bundle of a Riemannian manifold of positive
constant sectional curvature. The obtained class of K\"ahler
Einstein structures depends on one essential parameter, cannot
have constant holomorphic sectional curvature and is not
locally symmetric.}
\vskip3mm \scriptsize {\it Abbreviated title}:
A CLASS OF K\"AHLER EINSTEIN STRUCTURES
\end{flushleft}

\vskip5mm \centerline {\small \bf 1.INTRODUCTION} \vskip5mm

 In the study of the differential geometry of the cotangent bundle
$T^*M$ of a Riemannian manifold $(M,g)$ one uses several
Riemannian and semi-Riemannian metrics, induced from the
Riemannian metric $g$ on $M$. Next, one can get from $g$ some
natural almost complex structures on $T^*M$. The study of the
almost Hermitian structures induced from $g$ on $T^*M$ is an
interesting problem in the differential geometry of the cotangent
bundle.

 In [9] the authors have obtained a class of natural
K\"ahler Einstein structures $(G,J)$ of diagonal type induced on
$T^*M$ from the Riemannian metric $g$. The obtained K\"ahler
structures on $T^*M$ depend on two essential parameters $a_1$ and
$\lambda$, which are smooth functions depending on the energy
density $t$ on $T^*M$. In the case where the considered K\"ahler
structures are Einstein they get several situations in which  the
parameters $a_1,\lambda $ are related by some algebraic relations.
In the general case, $(T^*M,G,J)$ has constant holomorphic
curvature.

 In the present paper we study the singular case where the paramater
$a_1=A \sqrt{t}, ~A\in {\bf R}$. The class of the natural almost
complex structures $J$ on the nonzero cotangent bundle $T^*_0M$
that interchange the vertical and horizontal distributions depends
on one essential parameter $b_1$. This parameter is a smooth real
function depending on the energy density $t$ on $T^*_0M$. From the
integrability condition for $J$ it follows that the base manifold
$M$ must have constant curvature $c=\frac{A^2}{2}$ and the
parameter $b_1$ must fulfill the condition
$b_1>-\sqrt{\frac{c}{2t}}$.

 A class of natural Riemannian metrics $G$ of diagonal type on $T^*_0M$
is defined by four parameters $c_1, c_2, d_1, d_2$ which are
smooth functions of~$t$. From the condition for $G$ to be
Hermitian with respect to $J$ we get two sets of proportionality
relations, from which we can get the parameters $c_1, c_2, d_1,
d_2$ as functions depending on two new parameters $\lambda $ and
$\mu$.

 In the case where the fundamental $2$-form
$\phi$, associated to the class of complex structures $(G,J)$ is
closed, one finds that $\mu = \lambda^\prime$.

 Thus, we get
a class of K\"ahler structures $(G,J)$ on $T^*_0M$, depending on
one essential parameter $\lambda$.

 Finally, we find out the conditions under which  the obtained class
 of K\"ahler structures $(G,J)$ on $T^*_0M$ is  Einstein .

The manifolds, tensor fields and geometric objects we consider in
this paper, are assumed to be differentiable of class $C^{\infty}$
(i.e. smooth). We use the computations in local coordinates but
many results from this paper may be expressed in an invariant
form. The well known summation convention is used throughout this
paper, the range for the indices $h,i,j,k,l,r,s$ being
always${\{}1,...,n{\}}$ (see [3], [7], [9]). We shall denote by
${\Gamma}(T^*_0M)$ the module of smooth vector fields on $T^*_0M$.

\vskip5mm \centerline {\small \bf 2. SOME GEOMETRIC PROPERTIES OF
$T^*M$} \vskip5mm

Let $(M,g)$ be a smooth $n$-dimensional Riemannian manifold and
denote its cotangent bundle by $\pi :T^*M\longrightarrow M$.
Recall that there is a structure of a $2n$-dimensional smooth
manifold on $T^*M$, induced from the structure of smooth
$n$-dimensional manifold  of $M$. From every local chart
$(U,\varphi )=(U,x^1,\dots ,x^n)$  on $M$, it is induced a local
chart $(\pi^{-1}(U),\Phi )=(\pi^{-1}(U),q^1,\dots , q^n,$
$p_1,\dots ,p_n)$, on $T^*M$, as follows. For a cotangent vector
$p\in \pi^{-1}(U)\subset T^*M$, the first $n$ local coordinates
$q^1,\dots ,q^n$ are  the local coordinates $x^1,\dots ,x^n$ of
its base point $x=\pi (p)$ in the local chart $(U,\varphi )$ (in
fact we have $q^i=\pi ^* x^i=x^i\circ \pi, \ i=1,\dots n)$. The
last $n$ local coordinates $p_1,\dots ,p_n$ of $p\in \pi^{-1}(U)$
are the vector space coordinates of $p$ with respect to the
natural basis $(dx^1_{\pi(p)},\dots , dx^n_{\pi(p)})$, defined by
the local chart $(U,\varphi )$,\ i.e. $p=p_idx^i_{\pi(p)}$.

 An $M$-tensor field of type $(r,s)$ on $T^*M$ is defined by sets
of $n^{r+s}$ components (functions depending on $q^i$ and $p_i$),
with $r$ upper indices and $s$ lower indices, assigned to induced
local charts $(\pi^{-1}(U),\Phi )$ on $T^*M$, such that the local
coordinate change rule is that of the local coordinate components
of a tensor field of type $(r,s)$ on the base manifold $M$ (see
[5] for further details in the case of the tangent bundle). An
usual tensor field of type $(r,s)$ on $M$ may be thought of as an
$M$-tensor field of type $(r,s)$ on $T^*M$. If the considered
tensor field on $M$ is covariant only, the corresponding
$M$-tensor field on $T^*M$ may be identified with the induced
(pullback by $\pi $) tensor field on $T^*M$.

Some useful $M$-tensor fields on $T^*M$ may be obtained as
follows. Let $u,v:[0,\infty ) \longrightarrow {\bf R}$ be a smooth
functions and let $\|p\|^2=g^{-1}_{\pi(p)}(p,p)$ be the square of
the norm of the cotangent vector $p\in \pi^{-1}(U)$ ($g^{-1}$ is
the tensor field of type (2,0) having the components $(g^{kl}(x))$
which are the entries of the inverse of the matrix $(g_{ij}(x))$
defined by the components of $g$ in the local chart $(U,\varphi
)$). The components $u(\|p\|^2)g_{ij}(\pi(p))$, $p_i$,
$v(\|p\|^2)p_ip_j $ define $M$-tensor fields of types $(0,2)$,
$(0,1)$, $(0,2)$ on $T^*M$, respectively. Similarly, the
components $u(\|p\|^2)g^{kl}(\pi(p))$, $g^{0i}=p_hg^{hi}$,
$v(\|p\|^2)g^{0k}g^{0l}$ define $M$-tensor fields of type $(2,0)$,
$(1,0)$, $(2,0)$ on $T^*M$, respectively. Of course, all the
components considered above are in the induced local chart
$(\pi^{-1}(U),\Phi)$.

 The Levi Civita connection $\dot \nabla $ of $g$ defines a direct
sum decomposition

\begin{equation}
TT^*M=VT^*M\oplus HT^*M.
\end{equation}
of the tangent bundle to $T^*M$ into vertical distributions
$VT^*M= {\rm Ker}\ \pi _*$ and the horizontal distribution
$HT^*M$.

 If $(\pi^{-1}(U),\Phi)=(\pi^{-1}(U),q^1,\dots ,q^n,p_1,\dots ,p_n)$
is a local chart on $T^*M$, induced from the local chart
$(U,\varphi )= (U,x^1,\dots ,x^n)$, the local vector fields
$\frac{\partial}{\partial p_1}, \dots , \frac{\partial}{\partial
p_n}$ on $\pi^{-1}(U)$ define a local frame for $VT^*M$ over $\pi
^{-1}(U)$ and the local vector fields $\frac{\delta}{\delta
q^1},\dots ,\frac{\delta}{\delta q^n}$ define a local frame for
$HT^*M$ over $\pi^{-1}(U)$, where
$$
\frac{\delta}{\delta q^i}=\frac{\partial}{\partial
q^i}+\Gamma^0_{ih} \frac{\partial}{\partial p_h},\ \ \ \Gamma
^0_{ih}=p_k\Gamma ^k_{ih}
 $$
and $\Gamma ^k_{ih}(\pi(p))$ are the Christoffel symbols of $g$.

The set of vector fields $(\frac{\partial}{\partial p_1},\dots
,\frac{\partial}{\partial p_n}, \frac{\delta}{\delta q^1},\dots
,\frac{\delta}{\delta q^n})$ defines a local frame on $T^*M$,
adapted to the direct sum decomposition (1).

We consider
\begin{equation}
t=\frac{1}{2}\|p\|^2=\frac{1}{2}g^{-1}_{\pi(p)}(p,p)=\frac{1}{2}g^{ik}(x)p_ip_k,
\ \ \ p\in \pi^{-1}(U)
\end{equation}
the energy density defined by $g$ in the cotangent vector $p$. We
have $t\in [0,\infty)$ for all $p\in T^*M$.

From now on we shall work in a fixed local chart $(U,\varphi)$ on
$M$ and in the induced local chart $(\pi^{-1}(U),\Phi)$ on $T^*M$.

Now we shall present the following auxiliary result. \vskip5mm

{\footnotesize LEMMA 1}. \it If ~$n>1$ and $u,v$ are smooth
functions on $T^*M$ such that
$$
u g_{ij}+v p_ip_j=0,\ p\in \pi^{-1}(U)
$$
on the domain of any induced local chart on $T^*M$, then $u=0,\
v=0$.\rm \vskip 0.5cm

The proof is obtained easily by transvecting the given relation
with the components $g^{ij}$ of the tensor field $g^{-1}$ and
$g^{0j}$.\vskip5mm

{\it Remark}. From the relations of the type
$$
u g^{ij}+v g^{0i}g^{0j}=0,\ p\in \pi^{-1}(U),
$$
$$
u\delta ^i_j+vg^{0i} p_j=0,\ p\in \pi^{-1}(U),
$$
it is obtained, in a similar way, $u=v=0$ \rm.

\vskip5mm \centerline {\small \bf 3. A CLASS OF NATURAL COMPLEX
STRUCTURES}\centerline {\footnotesize \bf OF DIAGONAL TYPE ON
$T^*_0M$} \vskip5mm

The nonzero cotangent bundle $T^*_0M$ of Riemannian manifold
$(M,g)$ is defined by the formula: $T^*M$ minus zero section.
 Consider the real valued smooth functions $a_1,a_2,b_1,b_2$
defined on $(0,\infty)$. We define a class of natural almost
complex structures $J$ of diagonal type on $T^*_0M$ , expressed in
the adapted local frame by

\begin{equation}
J{\frac{\delta}{\delta
q^i}}=J^{(1)}_{ij}(p){\frac{\partial}{\partial p_j}},
~~~~~J{\frac{\partial}{\partial
p_i}}=-J_{(2)}^{ij}(p){\frac{\delta}{\delta q^j}}.
\end{equation}

where,
\begin{equation}
\left\{
\begin{array}{l}
J^{(1)}_{ij}(p)=a_1(t)g_{ij} + b_1(t)p_ip_j,
\\ \mbox{ } \\
J_{(2)}^{ij}(p)=a_2(t)g^{ij} + b_2(t)g^{0i}g^{0j},~~~A\in {\bf
R^*}.
\end{array}
\right.
\end{equation}

 In this paper we study the singular case where

\begin{equation}
a_1(t)=A \sqrt{t},~~~A\in {\bf R^*}.
\end{equation}

 The components $J^{(1)}_{ij},J_{(2)}^{ij}$ define symmetric
 $M$-tensor fields of types $(0,2),(2,0)$ on $T^*M$, respectively.

\vskip5mm {\footnotesize PROPOSITION 2}. \it The operator $J$
defines an almost complex structure on $T^*M$ if and only if

\begin{equation}
a_1a_2=1,~~~ (a_1+2tb_1)(a_2+2tb_2)=1.
\end{equation}
\vskip2mm

 Proof. \rm The relations are obtained easily from the property $J^2=-I$ of
$J$ and Lemma 1.

From the relations (5), (6) we can obtain the explicit expression
of the parameter $a_2, b_2$
\begin{equation}
a_2=\frac{1}{A \sqrt{t}},~~~~~~~~b_2=\frac{-b_1}{At(A+2
\sqrt{t}b_1)}.
\end{equation}

The obtained class of almost complex structures defined by the
tensor field $J$ on $T^*_0M$ is called \it class of natural almost
complex structures of diagonal type, \rm obtained from the
Riemannian metric $g$, by using the parameter $b_1$. We use the
word diagonal for these almost complex structures, since the
$2n\times 2n$-matrix associated to $J$, with respect to the
adapted local frame $(\frac{\delta}{\delta q^1},\dots
,\frac{\delta}{\delta q^n},\frac{\partial}{\partial p_1},\dots
,\frac{\partial}{\partial p_n})$ has two $n\times n$-blocks on the
second diagonal
\begin{displaymath}
J= \left(
\begin{array}{cc}
0 & -J_{(2)}^{ij} \\
J^{(1)}_{ij} & 0
\end{array}
\right).
\end{displaymath}

{\it Remark}. From the conditions (6) it follows that $a_1=A
\sqrt{t}$ and $a_2=\frac{1}{A\sqrt{t}}$ cannot vanish and have the
same sign. We assume that

\begin{equation}
A>0.
\end{equation}

 Similarly, from the conditions (6) it follows that $a_1+2tb_1$ and
$a_2+2tb_2$ cannot vanish and have the same sign.~We assume that
$a_1+2tb_1>0,~a_2+2tb_2>0~~\forall t>0$, i.e.

\begin{equation}
 A+2\sqrt{t}b_1>0~~~~~\forall t>0.
\end{equation}

Now we shall study the integrability of the class of natural
almost complex structures defined by $J$ on $T^*_0M$. To do this
we need the following well known formulas for the brackets of the
vector fields $\frac{\partial}{\partial p_i}, \frac{\delta}{\delta
q^i},~ i=1,...,n$
\begin{equation}
[\frac{\partial}{\partial p_i},\frac{\partial}{\partial
p_j}]=0,~~~[\frac{\partial}{\partial p_i},\frac{\delta}{\delta
q^j}]=\Gamma^i_{jk}\frac{\partial}{\partial p_k},~~~
[\frac{\delta}{\delta q^i},\frac{\delta}{\delta q^j}]
=R^0_{kij}\frac{\partial}{\partial p_k},
\end{equation}
where $R^h_{kij}(\pi(p))$ are the local coordinate components of
the curvature tensor field of $\dot \nabla$ on $M$ and
$R^0_{kij}(p)=p_hR^h_{kij}$ . Of course, the components
 $R^h_{kij}$, $R^0_{kij}$ define M-tensor fields of types
 (1,3), (0,3) on $T^*_0M$, respectively.\\

 Recall that the Nijenhuis tensor field $N$ defined by $J$ is given
by
$$
N(X,Y)=[JX,JY]-J[JX,Y]-J[X,JY]-[X,Y],\ \ \forall\ \ X,Y \in \Gamma
(T^*_0M).
$$
Then, we have $\frac{\delta}{\delta q^k}t =0,\
\frac{\partial}{\partial p_k}t = g^{0k}$. The expressions for the
components of $N$ can be obtained by a quite long, straightforward
computation, as follows\\

 {\footnotesize THEOREM 3. } {\it The Nijenhuis tensor field of the almost
complex structure $J$ on $T^*_0M$ is given by}
$$
\left\{
\begin{array}{l}
N(\frac{\delta}{\delta q^i},\frac{\delta}{\delta
q^j})=\{\frac{A^2}{2}
(\delta^h_ig_{jk}-\delta^h_jg_{ik})-R^h_{kij}\}p_h\frac{\partial}{\partial
p_k},
\\ \mbox{ } \\
N(\frac{\delta}{\delta q^i},\frac{\partial}{\partial
p_j})=J_{(2)}^{kl}J_{(2)}^{jr}\{\frac{A^2}{2}(\delta^h_ig_{rl}-
\delta^h_rg_{il})-R^h_{lir}\}p_h\frac{\delta}{\delta q^k},
\\ \mbox{ } \\
N(\frac{\partial}{\partial p_i},\frac{\partial}{\partial
p_j})=J_{(2)}^{ir}J_{(2)}^{jl}\{\frac{A^2}{2}(\delta^h_lg_{rk}-
\delta^h_rg_{lk})-R^h_{klr}\}p_h\frac{\partial}{\partial p_k}.
\end{array}
\right.
$$

\vskip2mm
 {\footnotesize THEOREM 4. } {\it The almost complex structure $J$
on $T^*_0M$ is integrable if and only if $(M,g)$ has positive
constant sectional curvature $c$ and
\begin{equation}
A=\sqrt{2c}.
\end{equation}

The function $b_1$ must fulfill the conditon
\begin{equation}
b_1>-\sqrt{\frac{c}{2t}}.
\end{equation}
\rm
 {\it Proof. }
  From the condition $N=0$, one obtains
$$
\{\frac{A^2}{2}
(\delta^h_ig_{jk}-\delta^h_jg_{ik})-R^h_{kij}\}p_h=0
$$

Differentiating with respect to $p_l$, it follows that the
curvature tensor field of $\dot \nabla$ has the expression
$$
R^l_{kij}=\frac{A^2}{2}(\delta^l_ig_{jk}- \delta^l_jg_{ik}).
$$
 Thus $(M,g)$ has positive constant sectional curvature
$c=\frac{A^2}{2}$. It follows that $A=\sqrt{2c}> 0$.

Conversely, if $(M,g)$ has positive constant sectional curvature
$c$ and $A$ is given by (11), one obtains
 in a straightforward way that $N = 0$.\\

Using by the relations (9),(11) we obtain the condition
(12).\\

The class of natural complex structures $J$ of diagonal type on
$T^*_0M$ depends on one essential parameter $b_1$. The components
of $J$ are given by
\begin{equation}
\left\{
\begin{array}{l}
J^{(1)}_{ij}=\sqrt{2ct}g_{ij}+b_1p_ip_j,
\\ \mbox{ } \\
J_{(2)}^{ij}=\frac{1}{\sqrt{2ct}} g^{ij}-
\frac{b_1}{2\sqrt{c}t(\sqrt{c}+ \sqrt{2t}b_1)}g^{0i}g^{0j}.
\end{array}
\right.
\end{equation}

\vskip5mm \centerline {\small \bf 4. A CLASS OF NATURAL HERMITIAN
STRUCTURES ON $T^*_0M$}\vskip5mm

Consider the following symmetric $M-$tensor fields on $T^*_0M$,
defined by the components
\begin{equation}
G^{(1)}_{ij}=c_1 g_{ij}+d_1p_ip_j,\ \ \
G_{(2)}^{ij}=c_2g^{ij}+d_2g^{0i} g^{0j},
\end{equation}
where $c_1,c_2,d_1,d_2$ are smooth functions depending on the
energy density $t\in (0,\infty)$.

Obviously, $G^{(1)}$ is of type $(0,2)$ and $G_{(2)}$ is of type
$(2,0)$. We shall assume that the matrices defined by $G^{(1)}$
and $G_{(2)}$ are positive definite. This happens if and only if
\begin{equation}
c_1>0,~c_2>0,~c_1+2td_1>0,~ c_2+2td_2>0~~\forall t>0.
\end{equation}
Then the following class of Riemannian metrics may be considered
on $T^*_0M$
\begin{equation}
G=G^{(1)}_{ij}dq^idq^j+G_{(2)}^{ij}Dp_iDp_j,
\end{equation}
where $Dp_i=dp_i-\Gamma^0_{ij}dq^j$ is the absolute (covariant)
differential of $p_i$ with respect to the Levi Civita connection
$\dot\nabla$ of $g$. Equivalently, we have
$$
G(\frac{\delta}{\delta q^i},\frac{\delta}{\delta
q^j})=G^{(1)}_{ij},~~G(\frac{\partial}{\partial
p_i},\frac{\partial}{\partial p_j})=G_{(2)}^{ij},~~
G(\frac{\partial}{\partial p_i},\frac{\delta}{\delta q^j})=
G(\frac{\delta}{\delta q^j},\frac{\partial}{\partial p_i})=0.
$$
Remark that $HT^*_0M,~VT^*_0M$ are orthogonal to each other with
respect to $G$, but the Riemannian metrics induced from $G$ on
$HT^*_0M,~VT^*_0M$ are not the same, so the considered metric $G$
on $T^*_0M$ is not a metric of Sasaki type. The $2n\times
2n$-matrix associated to $G$, with respect to the adapted local
frame $(\frac{\delta}{\delta q^1},\dots ,\frac{\delta}{\delta
q^n},\frac{\partial}{\partial p_1},\dots ,\frac{\partial}{\partial
p_n})$ has two $n\times n$-blocks on the first diagonal
\begin{displaymath}
G= \left(
\begin{array}{cc}
G^{(1)}_{ij} & 0  \\
0 & G_{(2)}^{ij}
\end{array}
\right).
\end{displaymath}

The class of Riemannian metrics $G$ is called a \it class of
natural lifts of diagonal type \rm of $g$.

 Remark also that the
system of 1-forms $(Dp_1,...,Dp_n,dq^1,...,dq^n)$ defines a local
frame on $T^{*}T^*_0M$, dual to the local frame
$(\frac{\partial}{\partial p_1},\dots ,\frac{\partial}{\partial
p_n}, \frac{\delta}{\delta q^1},\dots ,\frac{\delta}{\delta
q^n})$, on $TT^*_0M$ over $\pi^{-1}(U)$ adapted to the direct sum
decomposition (1).

We shall consider another two $M$-tensor fields $H_{(1)}, \
H^{(2)}$ on $T^*_0M$, defined by the components
$$
\left\{
\begin{array}{l}
H_{(1)}^{jk}=\frac{1}{c_1}g^{jk}-\frac{d_1}{c_1(c_1+2td_1)}g^{0j}g^{0k},
\\ \mbox{ } \\
H^{(2)}_{jk}=\frac{1}{c_2}g_{jk}-\frac{d_2}{c_2(c_2+2td_2)}p_jp_k.
\end{array}
\right.
$$

The components $H_{(1)}^{jk}$ define an $M$-tensor field of type
$(2,0)$ and the components $H^{(2)}_{jk}$ define an $M$-tensor
field of type $(0,2)$. Moreover, the matrices associated to
$H_{(1)}, \ H^{(2)}$ are the inverses of the matrices associated
to  $G^{(1)}$ and  $G_{(2)}$, respectively. Hence we have
$$
G^{(1)}_{ij}H_{(1)}^{jk} = \delta_i^k,\ \ G_{(2)}^{ij}H^{(2)}_{jk}
= \delta^i_k.
$$

Now, we shall be interested in the conditions under which the
class of the metrics $G$ is Hermitian with respect to the class of
the complex structures $J$, considered in the previous section,
i.e.
$$
G(JX,JY)=G(X,Y),
$$
for all vector fields $X,Y$ on $T^*_0M$.

Considering the coefficients of $g_{ij}, g^{ij}$ in the conditions
\begin{equation}
\left\{
\begin{array}{l}
G(J\frac{\delta}{\delta q^i},J\frac{\delta}{\delta q^j})=
G(\frac{\delta}{\delta q^i},\frac{\delta}{\delta q^j}),
\\ \mbox{ } \\
G(J\frac{\partial}{\partial p_i},J\frac{\partial}{\partial p_j})=
G(\frac{\partial}{\partial p_i},\frac{\partial}{\partial p_j}),
\end{array}
\right.
\end{equation}
we can express the parameters $c_1,c_2$  with the help of the
parameters $a_1, a_2$ and a proportionality factor $\lambda =
\lambda (t)$ . Then
\begin{equation}
c_1 = \lambda a_1 = \sqrt{2ct}\lambda,~~~~~~~~\ c_2=\lambda a_2 =
\frac{\lambda}{\sqrt{2ct}},
\end{equation}
where the coefficients $a_1,a_2$ are given by (5) and (7). Since
we made the assumption $c_1>0,\ c_2>0$, it follows $\lambda>0$.

Next, considering the coefficients of $p_ip_j,\ g^{0i}g^{0j}$ in
the relations (17), we can express the parameters
$c_1+2td_1,c_2+2td_2$ with help of the parameters $a_1+2tb_1,
a_2+2tb_2$ and a proportionality factor $\lambda+2t\mu $

\begin{equation}
\left\{
\begin{array}{l}
c_1+2td_1=(\lambda+2t\mu )(a_1+2tb_1),
\\ \mbox{ } \\
c_2+2td_2=(\lambda+2t\mu )(a_2+2tb_2).
\end{array}
\right.
\end{equation}
Remark that $\lambda (t)+2t \mu (t) > 0~ \forall t>0$. It is much
more convenient to consider the proportionality factor in such a
form in the expression of the parameters $c_1+2td_1,c_2+2td_2$.
Using by the relations (5), (7), (11),(18) we can obtain easily
from (19) the explicit expressions of the coefficients $d_1,d_2$
\begin{equation}
\left\{
\begin{array}{l}
d_1= b_1\lambda +  \sqrt{2t}(\sqrt{c} + \sqrt{2t} b_1)\mu,
\\ \mbox{ } \\
d_2=\frac{-b_1 \lambda +
\sqrt{2ct}\mu}{2\sqrt{c}t(\sqrt{c}+\sqrt{2t}b_1)}.
\end{array}
\right.
\end{equation}
Hence we may state:

{\footnotesize THEOREM 5.} \it Let $J$ be the class of natural,
complex structure of diagonal type on $T^*_0M$, given by (3) and
(13). Let $G$ be the class of the natural Riemannian metrics of
diagonal type on $T^*_0M$, given by (14), (18), (20).

Then we obtain a class of Hermitian structures $(G,J)$ on
$T^*_0M$, depending on three essential parameters $b_1$, $\lambda$
and $\mu$, which must fulfill  the conditions
\begin{equation}
b_1>-\sqrt{\frac{c}{2t}},~~~ \lambda > 0,~~~\lambda +2t\mu
>0~~~~\forall t>0.
\end{equation} \rm

\vskip5mm \centerline {\small \bf 5. A CLASS OF K\"AHLER
STRUCTURES ON $T^*_0M$} \vskip5mm
 Consider now the two-form $\phi $ defined by the class of
Hermitian structures $(G,J)$ on $T^*_0M$
$$
\phi (X,Y)=G(X,JY),
$$
for all vector fields $X,Y$ on $T^*_0M$.

Using by the expression of $\phi$ and computing the values $\phi
(\frac{\partial}{\partial p_i},\frac{\partial}{\partial p_j}),
\phi(\frac{\delta}{\delta q^i},\frac{\delta}{\delta q^j}),\\
\phi(\frac{\partial}{\partial p_i},\frac{\delta}{\delta q^j})$, we
obtain. \vskip5mm

 {\footnotesize PROPOSITION 6}. \it The expression of the
$2$-form $\phi $ in a local adapted frame
$(\frac{\partial}{\partial p_1},\dots ,\frac{\partial}{\partial
p_n}, \frac{\delta}{\delta q^1},\dots ,\frac{\delta}{\delta q^n})$
on $T^*_0M$, is given by
$$
\phi (\frac{\partial}{\partial p_i},\frac{\partial}{\partial
p_j})=0,\ \phi(\frac{\delta}{\delta q^i},\frac{\delta}{\delta
q^j})=0,\ \phi(\frac{\partial}{\partial p_i},\frac{\delta}{\delta
q^j})= \lambda \delta^i_j+\mu g^{0i}p_j,
$$
or, equivalently
\begin{equation}
\phi =(\lambda \delta^i_j+\mu g^{0i}p_j)Dp_i\wedge dq^j.
\end{equation}
\rm

{\footnotesize THEOREM 7}. \it The class of Hermitian structures
$(G,J)$ on $T^*_0M$ is K\"ahler if and only if
$$
\mu=\lambda ^\prime .
$$
Proof. \rm The expressions of $d\lambda,\ d\mu, \ dg^{0i}$ and
$dDp_i$ are obtained in a straightforward way, by using the
property $\dot \nabla_k g_{ij}=0$ (hence $\dot \nabla_k g^{ij}=0$)
$$
d\lambda = \lambda ^\prime g^{0i}Dp_i,\ d\mu =\mu ^\prime
g^{0i}Dp_i,\ dg^{0i}=g^{ik}Dp_k-g^{0h}\Gamma ^i_{hk}dq^k,
$$
$$
dDp_i=-\frac{1}{2}R^0_{ikl}dq^k\wedge dq^l+ \Gamma
^l_{ik}dq^k\wedge Dp_l.
$$
Then we have
$$
d\phi =(d\lambda \delta^i_j+d\mu g^{0i}p_j+ \mu dg^{0i}p_j+\mu
g^{0i}dp_j)\wedge Dp_i\wedge dq^j+
$$
$$
+(\lambda \delta^i_j+\mu g^{0i}p_j)dDp_i\wedge dq^j.
$$
By replacing the expressions of $d\lambda , d\mu , dg^{0i}$ and
$d\dot\nabla y^h$, then using, again, the property $\dot\nabla
_kg_{ij}=0$, doing some algebraic computations with the exterior
products, then using the well known symmetry properties of
$g_{ij}, \Gamma ^h_{ij},$ and of the Riemann-Christoffel tensor
field, as well as the Bianchi identities, it follows that
$$
d\phi =\frac{1}{2}(\lambda ^\prime -\mu)g^{0h}Dp_h\wedge
Dp_i\wedge dq^i.
$$
Therefore we have $d\phi =0$ if and only if $\mu =\lambda ^\prime
$. \vskip5mm

 {\it Remark}. The class of natural K\"ahler structures
of diagonal type defined by $(G,J)$ on $T^*_0M$ depends on two
essential parameters $b_1$ and $\lambda$.

The paramaters $b_1$ and $\lambda$ must fulfill the conditions
\begin{equation}
b_1>-\sqrt{\frac{c}{2t}},~~~\lambda > 0,~~~\lambda +
2t\lambda^{\prime}>0~~\forall t>0.
\end{equation}

 The components of the class of K\"ahler
metrics $G$ on $T^*_0M$ are given by
\begin{equation}
\left\{
\begin{array}{l}
G^{(1)}_{ij} =\sqrt{2ct}\lambda g_{ij}
+[\sqrt{2ct}\lambda^{\prime}+b_1(\lambda + 2t\lambda^{\prime})]
p_ip_j,
\\ \mbox{ } \\
G^{ij}_{(2)} =\frac{\lambda}{\sqrt{2ct}} g^{ij} +
\frac{-b_1\lambda+\sqrt{2ct}\lambda^{\prime}}{2\sqrt{c}t(\sqrt{c}+\sqrt{2t}b_1)}g^{0i}g^{0j}.
\end{array}
\right.
\end{equation}

We obtain, too
\begin{equation}
\left\{
\begin{array}{l}
H^{jk}_{(1)}=\frac{1}{\sqrt{2ct}\lambda}g^{jk} -
\frac{\sqrt{2ct}\lambda^{\prime}+b_1(\lambda +
2t\lambda^{\prime})}{2\sqrt{c}t\lambda(\sqrt{c}+\sqrt{2t}b_1)(\lambda
+ 2t\lambda^{\prime})}g^{0j}g^{0k},
\\ \mbox{ } \\
H^{(2)}_{jk}=\frac{\sqrt{2ct}}{\lambda}g_{jk}-
\frac{-b_1\lambda+\sqrt{2ct}\lambda^{\prime}}{\lambda(\lambda +
2t\lambda^{\prime})}p_jp_k.
\end{array}
\right.
\end{equation}
\vskip4cm

 \centerline {\small \bf 6. THE LEVI CIVITA CONNECTION OF
METRIC G}\centerline {\footnotesize \bf AND ITS CURVATURE TENSOR
FIELD} \vskip5mm

The Levi Civita connection $\nabla$ of the Riemannian manifold
$(T^*_0M,G)$ is determined by the conditions
$$
\nabla G=0,~~~~~  T =0,
$$
where $T$ is its torsion tensor field. The explicit expression of
this connection is obtained from the formula
$$
2G({\nabla}_XY,Z)=X(G(Y,Z))+Y(G(X,Z))-Z(G(X,Y))+
$$
$$
+G([X,Y],Z)-G([X,Z],Y)-G([Y,Z],X); ~~~~~~ \forall\
X,Y,Z~{\in}~{\Gamma}(T^*_0M).
$$

The final result can be stated as follows. \\

{\footnotesize THEOREM 8}. {\it The Levi Civita connection
${\nabla}$ of $G$ has the following expression in the local
adapted frame $(\frac{\delta}{\delta q^1},\dots
,\frac{\delta}{\delta q^n},\frac{\partial}{\partial p_1},\dots
,\frac{\partial}{\partial p_n}):$
\begin{equation}
\left\{
\begin{array}{l}
 \nabla_\frac{\partial}{\partial
p_i}\frac{\partial}{\partial p_j}
=Q^{ij}_h\frac{\partial}{\partial p_h},\ \ \ \ \ \
\nabla_\frac{\delta}{\delta q^i}\frac{\partial}{\partial
p_j}=-\Gamma^j_{ih}\frac{\partial}{\partial
p_h}+P^{hj}_i\frac{\delta}{\delta q^h},
\\ \mbox{ } \\
\nabla_\frac{\partial}{\partial p_i}\frac{\delta}{\delta
q^j}=P^{hi}_j\frac{\delta}{\delta q^h},\ \ \ \ \ \
\nabla_\frac{\delta}{\delta q^i}\frac{\delta}{\delta
q^j}=\Gamma^h_{ij}\frac{\delta}{\delta
q^h}+S_{hij}\frac{\partial}{\partial p_h},
\end{array}
\right.
\end{equation}

where $Q^{ij}_h, P^{hi}_j, S_{hij}$ are $M$-tensor fields on
$T^*_0M$, defined by
\begin{equation}
\left\{
\begin{array}{l}
Q^{ij}_h = \frac{1}{2}H^{(2)}_{hk}(\frac{\partial}{\partial
p_i}G_{(2)}^{jk}+ \frac{\partial}{\partial p_j}G_{(2)}^{ik}
-\frac{\partial}{\partial p_k}G_{(2)}^{ij}),
\\ \mbox{ } \\
P^{hi}_j=\frac{1}{2}H_{(1)}^{hk}(\frac{\partial}{\partial
p_i}G^{(1)}_{jk}-G_{(2)}^{il}R^0_{ljk}),
\\ \mbox{ } \\
S_{hij}=-\frac{1}{2}H^{(2)}_{hk}\frac{\partial}{\partial
p_k}G^{(1)}_{ij}+\frac{1}{2}R^0_{hij}.
\end{array}
\right.
\end{equation}

 \rm

Assuming that the base manifold $(M,g)$ has positive constant
sectional curvature $c$ and the class of metrics $G$ is K\"ahler,
using by the relations (24),(25), one obtains

$$
\left\{
\begin{array}{l}
 Q^{ij}_h =\frac{\sqrt{c}-\sqrt{2t}b_1}{4\sqrt{c}t}g^{ij}p_h-
\frac{\lambda-2t\lambda^{\prime}}{4t\lambda}(\delta^i_hg^{0j}+
\delta^j_hg^{0i})+
\\ \mbox{ } \\
~~~~~~~~\frac{\lambda(b_1^2\lambda+c\lambda^{\prime})-\sqrt{2c}\lambda(b_1^{\prime}\lambda-b_1\lambda^{\prime})t^{1/2}
+2(b_1^{\prime}\lambda\lambda^{\prime}-2c\lambda^{\prime2}+c\lambda\lambda^{\prime\prime})t
-2\sqrt{2c}\lambda^{\prime}(\lambda
b_1^{\prime}+2b_1\lambda^{\prime})t^{3/2}}{4\sqrt{c}t\lambda(\sqrt{c}+\sqrt{2t}b_1)(\lambda+2t\lambda^{\prime})}g^{0i}g^{0j}p_h,
\\ \mbox{ } \\
P^{hi}_j=\frac{-\sqrt{c}\lambda+\sqrt{2}b_1\lambda
t^{1/2}+2\sqrt{c}\lambda^{\prime}t+2\sqrt{2}b_1\lambda^{\prime}t^{3/2}}{4\sqrt{c}\lambda
t}g^{hi}p_j+
\\ \mbox{ } \\
~~~~~~~~\frac{1}{4t}\delta^i_j
g^{0h}+\frac{\lambda+2t\lambda^{\prime}}{4t\lambda}\delta^h_jg^{0i}+\frac{1}{2\sqrt{2c}\lambda
t(\sqrt{2c}+2b_1\sqrt{t})(\lambda+2t\lambda^{\prime})}[\lambda(-b_1^2\lambda+c\lambda^{\prime})+
\\ \mbox{ } \\
~~~~~~~~\sqrt{2c}b_1^{\prime}\lambda^2
t^{1/2}-2(2b_1^2\lambda\lambda^{\prime}+2c\lambda^{\prime2}-
c\lambda\lambda^{\prime\prime})t+
\\ \mbox{ } \\
~~~~~~~~2\sqrt{2c}(b_1^{\prime}\lambda\lambda^{\prime}-3b_1\lambda^{\prime2}+
b_1\lambda\lambda^{\prime\prime})t^{3/2}-4b_1^2\lambda^{\prime2}t^2]
g^{0h}g^{0i}p_j,
\\ \mbox{ } \\
S_{hij}=\frac{\sqrt{c}(\sqrt{c}+\sqrt{2t}b_1)}{2}g_{ij}p_h
+\frac{\sqrt{c}(\sqrt{c}\lambda-\sqrt{2}b_1\lambda
t^{1/2}-2\sqrt{c}\lambda^{\prime}t-2\sqrt{2}b_1\lambda^{\prime}t^{3/2})}{2\lambda}g_{hj}p_i-
\\ \mbox{ } \\
~~~~~~~~\frac{\sqrt{c}(\sqrt{c}\lambda+\sqrt{2}b_1\lambda
t^{1/2}+2\sqrt{c}\lambda^{\prime}t+2\sqrt{2}b_1\lambda^{\prime}t^{3/2})}{2\lambda}g_{hi}p_j-
\\ \mbox{ } \\
~~~~~~~~\frac{1}{2\lambda(\lambda+2t\lambda^{\prime})}[\lambda(2b_1^2\lambda+c\lambda^{\prime})+
\sqrt{2c}\lambda(\lambda
b_1^{\prime}+4\lambda^{\prime}b_1)t^{1/2}+
\\ \mbox{ } \\
~~~~~~~~2(\lambda^2b_1b_1^{\prime}+5\lambda
\lambda^{\prime}b_1^2-2c\lambda^{\prime2}+c\lambda
\lambda^{\prime\prime})t+
\\ \mbox{ } \\
~~~~~~~~2\sqrt{2c}(\lambda
\lambda^{\prime}b_1^{\prime}-2\lambda^{\prime2}b_1+2\lambda
\lambda^{\prime\prime}b_1)t^{3/2}+ 4\lambda
b_1(\lambda^{\prime}b_1^{\prime}+\lambda^{\prime\prime}b_1)t^2]
p_hp_ip_j.
\end{array}
\right.
$$

 The curvature tensor field $K$ of the connection
$\nabla $ is obtained from the well known formula
$$
K(X,Y)Z=\nabla_X\nabla_YZ-\nabla_Y\nabla_XZ-\nabla_{[X,Y]}Z,\ \ \
\ \forall\ X,Y,Z\in \Gamma (T^*_0M).
$$

The components of curvature tensor field $K$ with respect to the
adapted local frame $(\frac{\delta}{\delta q^1},\dots
,\frac{\delta}{\delta q^n},\frac{\partial}{\partial p_1},\dots
,\frac{\partial}{\partial p_n})$ are obtained easily:
\begin{equation}
\left\{
\begin{array}{l}
K(\frac{\delta}{\delta q^i},\frac{\delta}{\delta
q^j})\frac{\delta}{\delta q^k}=QQQ^h_{ijk}\frac{\delta}{\delta
q^h},\ \ \ K(\frac{\delta}{\delta q^i},\frac{\delta}{\delta
q^j})\frac{\partial}{\partial
p_k}=QQP^k_{ijh}\frac{\partial}{\partial p_h},
\\ \mbox{ } \\
K(\frac{\partial}{\partial p_i},\frac{\partial}{\partial
p_j})\frac{\delta}{\delta q^k}=PPQ^{ijh}_k\frac{\delta}{\delta
q^h},\ \ \ K(\frac{\partial}{\partial
p_i},\frac{\partial}{\partial p_j})\frac{\partial}{\partial p_k}
=PPP^{ijk}_h\frac{\partial}{\partial p_h},
\\ \mbox{ } \\
K(\frac{\partial}{\partial p_i},\frac{\delta}{\delta
q^j})\frac{\delta}{\delta q^k}=PQQ^i_{jkh}\frac{\partial}{\partial
p_h},\ \ \ K(\frac{\partial}{\partial p_i},\frac{\delta}{\delta
q^j})\frac{\partial}{\partial p_k}=PQP^{ikh}_j\frac{\delta}{\delta
q^h},
\end{array}
\right.
\end{equation}
where
\begin{equation}
\left\{
\begin{array}{l}
QQQ^h_{ijk}=R^h_{kij}-P^{hl}_kR^0_{lij}+P^{hl}_iS_{ljk}-P^{hl}_jS_{lik},
\\ \mbox{ } \\
QQP^k_{ijh}=-R^k_{hij} - Q^{lk}_h
R^0_{lij}-P^{lk}_iS_{hjl}+P^{lk}_j S_{hil},
\\ \mbox{ } \\
PPQ^{ijh}_k=\frac{\partial}{\partial
p_i}P^{hj}_k-\frac{\partial}{\partial
p_j}P^{hi}_k+P^{hi}_lP^{lj}_k-P^{hj}_lP^{li}_k,
\\ \mbox{ } \\
PPP^{ijk}_h=\frac{\partial}{\partial
p_i}Q^{jk}_h-\frac{\partial}{\partial
p_j}Q^{ik}_h+Q^{jk}_lQ^{il}_h- Q^{ik}_lQ^{jl}_h,
\\ \mbox{ } \\
PQQ^i_{jkh}=\frac{\partial}{\partial
p_i}S_{hjk}+Q^{il}_hS_{ljk}-P^{li}_kS_{hjl},
\\ \mbox{ } \\
PQP^{ikh}_j=\frac{\partial}{\partial
p_i}P^{hk}_j+P^{hi}_lP^{lk}_j-P^{hl}_jQ^{ik}_l.
\end{array}
\right.
\end{equation}
are M-tensor fields on $T^*_0M$.\\

 {\it Remark.} The explicit expressions of these components are obtained after
some quite long and hard computations, made by using the package
RICCI.

From the local coordinates expression of the curvature tensor
field K we obtain that the class of  K\"ahler structures $(G,J)$
on $T^*_0M$ cannot have constant holomorphic
sectional curvature.\\

The Ricci tensor field Ric of $\nabla$ is defined by the formula:
$$
Ric(Y,Z)=trace(X\longrightarrow K(X,Y)Z),\ \ \ \forall\  X,Y,Z\in
\Gamma (T^*_0M).
$$
It follows
$$
\left\{
\begin{array}{l}
Ric(\frac{\delta}{\delta q^j},\frac{\delta}{\delta
q^k})=RicQQ_{jk}=QQQ_{hjk}^h+PQQ_{jkh}^h,
\\ \mbox{ } \\
 Ric(\frac{\partial}{\partial
p_j},\frac{\partial}{\partial
p_k})=RicPP^{jk}=PPP^{hjk}_h-PQP^{jkh}_h,
\\ \mbox{ } \\
Ric(\frac{\partial}{\partial p_i},\frac{\delta}{\delta
q^j})=Ric(\frac{\delta}{\delta q^j},\frac{\partial}{\partial
p_i})=0.
\end{array}
\right.
$$

Doing the necessary computations, we obtain the final expressions
of the components of the Ricci tensor field of $\nabla$

\begin{equation}
\left\{
\begin{array}{l}
RicQQ_{jk}=\frac{a}{2\lambda(\lambda+2t\lambda^{\prime})}g_{jk}+\frac{\alpha}{4t\lambda^2
(\lambda+2t\lambda^{\prime})^2}p_jp_k,
\\ \mbox{ } \\
RicPP^{jk}=\frac{a}{4ct\lambda(\lambda+2t\lambda^{\prime})}g^{jk}+
\frac{\beta}{8\sqrt{c}t^2\lambda^2(\sqrt{c}+\sqrt{2t}b_1)(\lambda+2t\lambda^{\prime})^2}g^{0j}g^{0k},
\end{array}
\right.
\end{equation}
where the coefficient $a$ is given by
$$
a=c\lambda^2(n-2)-8c\lambda \lambda^{\prime}t
-4c(n-1)\lambda^{\prime2}t^2-4c\lambda \lambda^{\prime\prime}t^2-
$$
$$
\sqrt{2ct}[(n+1)\lambda^2+2\lambda
\lambda^{\prime}(2n+3)t+4\lambda^{\prime2}(n-1)t^2+4\lambda
\lambda^{\prime\prime}t^2]b_1-
$$
$$
2\sqrt{2c}\lambda t^{3/2}(\lambda+2t\lambda^{\prime})b_1^{\prime}.
$$

The expressions of $\alpha$ and $\beta$ are too large to be
presented in this paper.

\vskip5mm \centerline {\small \bf 7. A CLASS OF K\"AHLER EINSTEIN
STRUCTURES ON $T^*_0M$} \vskip5mm

In order to find out the conditions under which the class of
K\"ahler structures $(G,J)$ on $T^*_0M$ is Einstein, we consider
the relations

\begin{equation}
\left\{
\begin{array}{l}
RicQQ_{jk}=Ef~G^{(1)}_{jk},
\\ \mbox{ } \\
RicPP^{jk}=Ef~G_{(2)}^{jk},
\end{array}
\right.
\end{equation}
where the real number $Ef$ is an Einstein factor.

Using by the relations (24), (30) and Lemma 1, one obtains the
conditions

$$
Ef=\frac{a}{2\lambda^2\sqrt{2ct}(\lambda+2t\lambda^{\prime})}=
\frac{\alpha}{4t\lambda^2(\lambda+2t\lambda^{\prime})^2[\sqrt{2ct}\lambda^{\prime}+(\lambda+2t\lambda^{\prime})b_1]}=
$$
$$
=\frac{\beta}{4t\lambda^2(\lambda+2t\lambda^{\prime})^2(\sqrt{2ct}\lambda^{\prime}-\lambda
b_1)}.
$$

The condition
$$
Ef=\frac{a}{2\lambda^2\sqrt{2ct}(\lambda+2t\lambda^{\prime})}
$$
is equivalent with the first order differential equation
$$
2\sqrt{2c}\lambda t^{3/2}(\lambda+2t\lambda^{\prime})b_1^{\prime}=
-\sqrt{2ct}[(n+1)\lambda^2+2\lambda
\lambda^{\prime}(2n+3)t+4\lambda^{\prime2}(n-1)t^2+4\lambda
\lambda^{\prime\prime}t^2]b_1+
$$
$$
c\lambda^2(n-2)-8c\lambda \lambda^{\prime}t
-4c(n-1)\lambda^{\prime2}t^2-4c\lambda \lambda^{\prime\prime}t^2
-Ef~2\sqrt{2ct}\lambda^2(\lambda+2t\lambda^{\prime}).
$$

The solution of this equation is
\begin{equation}
b_1(t)=\frac{t^{-\frac{n+1}{2}}\lambda^{1-n}}{4(\lambda+2t\lambda^{\prime})}(4C-\int_{1}^{t}\theta(s)ds),
~~~t>0,~~~C\in {\bf R},
\end{equation}
where
$$
\theta(s)=s^{\frac{n-2}{2}}\lambda^{n-2} \{4~Ef
\sqrt{s}\lambda^2(\lambda+2s\lambda^{\prime})+\sqrt{2c}[-(n-2)\lambda^2
+4(n-1)s^2\lambda^{\prime2}+4s\lambda(2\lambda^{\prime}+s\lambda^{\prime\prime})]
\}.
$$
If the function $b_1$ is given by (32) then the relations
$$
Ef=
\frac{\alpha}{4t\lambda^2(\lambda+2t\lambda^{\prime})^2[\sqrt{2ct}\lambda^{\prime}+(\lambda+2t\lambda^{\prime})b_1]}=
\frac{\beta}{4t\lambda^2(\lambda+2t\lambda^{\prime})^2(\sqrt{2ct}\lambda^{\prime}-\lambda
b_1)}
$$
are fulfilled.

Now we may state our main result. \vskip5mm
 {\footnotesize THEOREM 9.}
{\it Assume that the Riemannian manifold $(M,g)$ has positive
constant sectional curvature $c$. Let $J$ be the class of natural,
complex structure of diagonal type on $T^*_0M$, given by (3) and
(13). Let $G$ be the class of the natural Riemannian metrics of
diagonal type on $T^*_0M$, given by (14) and (24).

If ~$b_1$ is given by (32) then $(G,J)$ is a class of K\"ahler
Einstein structures on $T^*_0M$, depending on one essential
parameter $\lambda$.

The function $b_1$ and the parameter $\lambda$ must fulfill the
conditions (23):
$$
b_1>-\sqrt{\frac{c}{2t}},~~~~\lambda > 0,~~~~\lambda +
2t\lambda^{\prime}>0~~~~\forall t>0.
$$
\rm

{\it Remark.} After some long and hard computations we have
obtained that the class of K\"ahler Einstein structures $(G,J)$ on
$T^*_0M$ is not locally symmetric.

{\it Example.} The case where $\lambda = 1$ is presented in [10].

\vskip1cm
 \centerline{\small REFERENCES}
\vskip1cm

\begin{flushleft}
[1]~A.~Besse, {\it Einstein manifolds.}~Ergeb.Math.Grenzgeb.(3)
10, Springer-Verlag, Berlin 1987.\vskip4mm

[2]~E. Calabi, {\it M\'etriques Kaehl\'eriennes et fibr\'es
holomorphes.}~Ann. Scient. Ec. Norm. Sup.,~12
(1979),~269-294.\vskip4mm

[3]~Gh. Gheorghiev and V. Oproiu, {\it Variet\u a\c ti diferen\c
tiabile finit \c si infinit dimensionale.}~Ed. Academiei Rom. I
(1976),~~II (1979).\vskip4mm

[4]~J.M. Lee, {\it Ricci. A Mathematica package for doing tensor
calculations in differential geometry. User's Manual. 1992,
2000.}\vskip4mm

[5]~K.P. Mok, E.M. Patterson and Y.C. Wong, {\it Structure of
symmetric tensors of type (0,2) and tensors of type (1,1) on the
tangent bundle}.~ Trans. Am. Math. Soc.~234
(1977),~253-278.\vskip4mm

[6]~V. Oproiu, {\it Some new geometric structures on the tangent
bundle.} Public. Math. Debrecen, 55 (1999), 261-281.\vskip4mm

[7]~V. Oproiu and N. Papaghiuc, {\it A K\" ahler structure on the
nonzero tangent bundle a space form.} Differential Geom. Appl. 11
(1999), 1-12.\vskip4mm

[8]~V. Oproiu and D.D. Poro\c sniuc, {\it A K\"ahler Einstein
structure on the cotangent bundle of a Riemannian manifold.}~An.
\c Stiin\c t. Univ. Al. I. Cuza, Ia\c si, 49 (2003), f2,
399-414.\vskip4mm

[9]~V. Oproiu and D.D. Poro\c sniuc, {\it A class of K\"ahler
Einstein structures on the cotangent bundle of a space form.}
Public. Math. Debrecen. To appear.\vskip4mm

[10]~D.D. Poro\c sniuc, {\it A K\" ahler Einstein structure on the
nonzero cotangent bundle of a space form.} Italian Journal of Pure
and Applied Mathematics. To appear.\vskip4mm

[11]~K. Yano and S. Ishihara,
{\it Tangent and Cotangent Bundles.}~M. Dekker Inc., New%
York,~1973.\vskip3cm

\end{flushleft}

\begin{flushright}
\scriptsize \it Department of Mathematics \\National College "M. Eminescu" \\
Str. Octav Onicescu 52 \\RO-710096 Boto\c sani, Rom\^ ania.\\
e-mail: dporosniuc@yahoo.com \\
~~~~~~~~~~danielporosniuc@lme.ro
\end{flushright}

\end{document}